\documentclass{fic-l}

\newtheorem{theorem}{Theorem}[section]

\theoremstyle{definition}
\newtheorem{definition}[theorem]{Definition}

\newtheorem{conjecture}[theorem]{Conjecture}
\newtheorem{problem}[theorem]{Problem}

\theoremstyle{remark}

\numberwithin{equation}{section}

\usepackage{amscd}

\usepackage{amssymb}
\usepackage{pb-diagram} 
\usepackage{graphicx} 

\newcommand{\TryPackage}[3]{\IfFileExists{#1.sty}{\usepackage{#1}#2}{#3}}
\TryPackage{mathrsfs}{\renewcommand{\mathcal}{\mathscr}}{%
     \TryPackage{eucal}{}{}}

\newcommand{\al}{\alpha}
\newcommand{\be}{\beta}

\newcommand{\varep}{\varepsilon}

\newcommand{\la}{\lambda}

\newcommand{\si}{\sigma}

\newcommand{\Si}{\Sigma}

\newcommand{\ZZ}{{\mathbb Z}}

\newcommand{\CC}{{\mathbb C}}

\newcommand{\SLC}{{SL_2({\mathbb C})}}
\newcommand{\PSLC}{{PSL_2({\mathbb C})}}

\newcommand{\lcm}{\operatorname{lcm}}

\begin{document}

\title{A { PSL}$_2(\CC)$ Casson invariant}

\author{Cynthia L. Curtis}
\address{Department of Mathematics and Statistics \\The College of New Jersey \\Ewing, NJ 08628}
\email{ccurtis@tcnj.edu}

\date{\today}

\begin{abstract}
We use intersection theory techniques to define an invariant of closed 3-manifolds counting the characters of irreducible representations of the fundamental group in $\PSLC$. We note several properties of the invariant and compute the invariant for certain Seifert fibered spaces and for some Dehn surgeries on twist knots.  We discuss the relationship between this invariant and the $\SLC$-invariant defined in \cite{C}.
\end{abstract}

\maketitle

\section{Introduction}

The purpose of this paper is to use intersection theory techniques to define an invariant of closed 3-manifolds which will count the equivalence classes of irreducible representations of the fundamental group in $\PSLC$.  The invariant defined here may be viewed as a $\PSLC$ analog of the Casson invariant.

Recall that the Casson invariant counts equivalence classes of irreducible representations of the fundamental group in $SU(2)$.  This invariant has been generalized to count representations in a number of other groups.  Of particular interest for us are the generalizations to $SO(3)$ found in \cite{C2} and to $\SLC$ found in \cite{C}. The $SO(3)$ invariant, like the Casson invariant, has a surgery formula for rational homology spheres, as shown in \cite{C2}.  A surgery formula for most surgeries on knots in homology spheres for the $\SLC$-invariant can be found in \cite{C}, as corrected in \cite{C1}, and further calculations for the $\SLC$-invariant can be found in \cite{BC}.

 Recall that $SO(3)$ is isomorphic to $PSU(2) = SU(2)/\pm I$.  In \cite{C2}, this fact is used to relate the piece of the $SO(3)$-invariant which counts representations which lift to $SU(2)$ to the $SU(2)$-invariant, and eventually to compute the $SO(3)$-invariant.  More specifically, we define the $SO(3)$-invariant of a closed 3-manifold $\Si$ to be a sum of invariants $\la_{\al}(\Si)$, where $\al \in H^2(\Si; \ZZ_2)$ and $\la_{\al}(\Si)$ counts the equivalence classes of representations $\rho$ with $w_2(\rho) = \al$. It is easy to check that $\la_0(\Si) = \frac{1}{|H^1(\Si;\ZZ_2)|}\la_{SU(2)}(\Si)$.  In fact, in the end it is shown that $\la_{\al}(\Si) = \la_0 (\Si)$ for every $\al \in H^2(\Si; \ZZ_2)$.  Since $H^1(\Si; \ZZ_2)$ and $H^2(\Si; \ZZ_2)$ are isomorphic, it follows that the $SO(3)$-invariant is equal to the the $SU(2)$-invariant.

Here, we will relate the portion of the $\PSLC$-invariant which counts representations which lift to $\SLC$ to the invariant of \cite{C} in a similar fashion.  Again, we may define invariants $\la_{\al}(\Si)$ counting the equivalence classes of irreducible $\PSLC$-representations $\rho$ with $w_2(\rho) = \al$. As in the case of $SO(3)$, we will see that $\la_0 (\Si) = \frac{1}{|H^1(\Si;\ZZ_2)|}\la_{\SLC)}(\Si)$. Interestingly, however, in contrast to the $SO(3)$-theory, we will show that in general $\la_{\al}(\Si) \neq \la_0 (\Si)$, so the $\PSLC$-invariant does {\em not} equal the $\SLC$-invariant.

In Section 2, we define the invariant $\la_\PSLC(\Si)$, and in Section 3 we establish several of its properties.  In Section 4 we compute the $\PSLC$- and $\SLC$-invariants for certain Seifert fibered spaces and for certain surgeries on twist knots, and in Section 5 we state several observations, questions, and conjectures concerning the invariant and its relationship to the $\SLC$-invariant.

\section{Definition of the $\PSLC$-invariant}

Given a finitely generated group $\pi$, denote by $ R(\pi)$ the
space of representations $\rho\colon \pi \to \PSLC$. It is well known
 that $R(\pi)$ has the structure of a complex affine
algebraic set.  (See \cite{BZ1}, for example.) Let $ X(\pi)$ denote the algebro-geometric quotient under the natural action of $\PSLC$ on $ R(\pi)$. Note that there is a surjective quotient map $ t: R(\pi) \to  X(\pi)$.  By analogy with the $\SLC$-theory, we call the elements of $ X(\pi)$ the $\PSLC$-characters of $\pi$.

 A representation $\rho \in  R(\pi)$ is {\em reducible} if it is conjugate to a representation whose image
 is upper triangular; all other representations are {\em irreducible}.  Let
$ R^*(\pi)$ denote the subspace of $R(\pi)$ consisting of irreducible representations. Let  $ X^*(\pi)$ be the subspace of $X(\pi)$ of characters of
irreducible representations. Given a manifold $\Si$, we denote by
$R(\Si)$ the variety of $\PSLC$ representations of $\pi_1 \Si$ and
by $ X(\Si)$ the associated character variety.

Suppose now $\Si$ is a closed, orientable 3-manifold
with a Heegaard splitting $(W_1, W_2, F)$. Here, $F$
is a closed orientable surface embedded in $\Si$, and $W_1$ and $W_2$ are handlebodies
with boundaries $\partial W_1 = F = \partial W_2$
 such that $\Si = W_1\cup_F W_2$.
The inclusion maps $F \hookrightarrow W_i$ and $W_i \hookrightarrow \Si$
induce a pushout diagram
$$\begin{diagram} 
\node{} \node{\pi_1 W_1}\arrow{se,l}{}\\
\node{\pi_1 F}\arrow{ne,l}{} \arrow{se,l}{} \node{} \node{\pi_1 \Si.}\\
\node{} \node{\pi_1 W_2} \arrow{ne,l}{}
\end{diagram}
$$
This diagram induces a diagram of the associated representation
varieties, where all arrows are reversed and are injective rather
than surjective. On the level of character varieties, this gives a
pullback diagram
 $$\begin{diagram} 
\node{} \node{X(W_1)}\arrow{sw,l}{}\\
\node{X(F)} \node{} \node{X(\Si).}\arrow{nw,l}{} \arrow{sw,l}{} \\
\node{} \node{X(W_2)} \arrow{nw,l}{}
\end{diagram}
$$
This identifies $X(\Si)$  as the intersection
$$X(\Si) = X(W_1) \cap X(W_2) \subset X(F).$$

To generalize the Casson invariant to $\PSLC$, we should define an intersection number for $X(W_1)$ and $X(W_2)$ in $X(F)$.  However the spaces involved are not manifolds, so we must do this carefully.  To begin with, the spaces are stratified according to the reducibility of the characters.  Since we are primarily interested in counting irreducible characters, we may avoid this problem by defining an intersection number which counts only those characters which lie in compact subvarieties of $X^*(F)$.  This is completely analogous to the approach taken in \cite{C} to define an invariant counting irreducible $\SLC$ characters.

However, unlike the $\SLC$ theory the variety $X^*(F)$ is not a manifold, but rather an orbifold.  Therefore we must define an orbifold intersection number, counting points in the intersection with weight $\pm 1/n$, where $n$ is the order of the stabilizer group.  This is completely analogous to the approach taken in \cite{C2} to define an invariant counting equivalence classes of irreducible $SO(3)$-representations.

We now describe the orbifold structure on $X^*(F)$.
Let $D \subset \PSLC$ denote the subgroup consisting of the images of the diagonal matrices and the
matrices $\pm \left( \begin{array}{cc}
    0 & z \\
    -z^{-1} & 0
    \end{array}\right)$
    where $z \in \Bbb C$, $z \neq 0$.
Let $Z$ denote the $\ZZ_2 \oplus \ZZ_2$-subgroup of $D$. For an irreducible representation $\rho$ which is conjugate to a representation into $D - Z$, it is easy to check that
$\PSLC$ does not act freely on the orbit $O(\rho)$, but rather acts with $\ZZ_2$-isotropy. Similarly, if $\rho$ is conjugate to a representation onto $Z$, then $\PSLC$  acts with $\ZZ_2 \oplus \ZZ_2$-isotropy. However,
$\PSLC$ does act freely on the orbits of all other irreducible representations.

There are natural orientations on all the character varieties
determined by their complex structures. The  invariant
$\la_\PSLC(\Si)$ is defined as an oriented orbifold intersection number of
$X^*(W_1)$ and $X^*(W_2)$ in $X^*(F)$ which counts only compact,
zero-dimensional components of the intersection.

Specifically,
there exist a compact neighborhood $U$ of the zero-dimensional
components of $X^*(W_1)\cap X^*(W_2)$ which is disjoint from the
higher dimensional components of the intersection and an isotopy
$h\colon X^*(F) \to X^*(F)$ supported in $U$ such that
$h(X^*(W_1))$ and $X^*(W_2)$ intersect transversely in $U$. This isotopy may not preserve the lower strata consisting of characters of representations into $D$ and into $Z$.  To account for this, we count a character $\chi$ in $h(X^*(W_1))\cap X^*(W_2)$ with weight $\be_{\chi} =1$ if $\chi$ is not the character of a representation that conjugates into $D$, with weight $\be_{\chi} = 1/2$ if $\chi$ is the character of a representation which conjugates into $D$ but not into $Z$, and with weight $\be_{\chi}=1/4$ if $\chi$ is the character of a representation which conjugates into $Z$.  Now set $\varep_\chi = \pm \be_\chi$, depending on whether the orientation of $h(X^*(W_1))$  followed by
that of $X^*(W_2)$ agrees with or disagrees with the orientation
of $X^*(F)$ at $\chi$.

\begin{definition} Let
$\la_\PSLC(\Si) = \sum_\chi \varep_\chi,$
where the sum is  over all zero-dimensional
components $\chi$ of the intersection $h(X^*(W_1))\cap X^*(W_2)$.
\end{definition}

Standard arguments show that the quantity $\la_\PSLC(\Si)$ is independent of the choice of isotopy used to define it
and of the choice of Heegaard splitting of $\Si$.  (See \cite{C} and \cite{C2}, for example.) Thus we obtain
\begin{theorem}
$\la_\PSLC(\Si)$ is a well-defined invariant of $\Si$.
\end{theorem}

Also as in the $\SLC$-theory, this definition agrees with the algebro-geometric intersection number.
Thus, for $\chi$ a zero-dimensional component of the intersection $X^*(W_1) \cap X^*(W_2) \subset X^*(F)$, if we define $n_{\chi}$ to be the coefficient of $\chi$ in the intersection cycle $X^*(W_1)\cdot X^*(W_2)$, we obtain
\begin{theorem}
$\la_\PSLC(\Si) = \sum_{\chi} n_{\chi},$ where the sum is taken over all zero-dimensional components of $X^*(W_1) \cap X^*(W_2)$.
\end{theorem}
\noindent See \cite{C} for details.

\section{Properties of the invariant}

We now establish some basic properties of the
$\PSLC$ Casson invariant. These follow quite easily from the definition and the relationship between $\SLC$ and $\PSLC$.
We show
\begin{theorem}
The invariant $\la_\PSLC$ has the following properties:
\begin{enumerate}
\item[(i)] For any 3-manifold $\Si,$ $\la_\PSLC(\Si)\geq 0.$  If $X^*(W_1) \cap X^*(W_2)$ has at least one zero-dimensional component, then
$\la_\PSLC(\Si) > 0.$
\item[(ii)] If $\Si$ is hyperbolic, then $\la_\PSLC(\Si) >0.$
\item[(iii)] If $\la_\PSLC(\Si) > 0,$ then there exists an
irreducible representation $\rho\colon \pi_1 \Si \to \PSLC$. In
particular, if $\pi_1 \Si$ is cyclic, then $\la_\PSLC(\Si) = 0.$
\item[(iv)] If $\Si$ is a $\ZZ_2$ homology 3-sphere, then $\la_\PSLC(\Si) = \la_\SLC(\Si).$
\item[(v)]  $\la_\PSLC(-\Si) =
\la_\PSLC(\Si)$, where $-\Si$ is $\Si$ with the opposite
orientation.
\item[(vi)] For rational homology spheres $\Si_1$ and $\Si_2$, we have
$$\la_\PSLC(\Si_1 \# \Si_2) =\la_\PSLC(\Si_1) + \la_\PSLC(\Si_2).$$
\end{enumerate}
\end{theorem}

\begin{proof}

\begin{enumerate}
\item[(i)] This is an immediate consequence of Definition 2.1 and Theorem 2.3.
\item[(ii)] The hyperbolic structure on $\Si$ is a discrete, faithful, irreducible representation $\pi_1(\Si) \to \PSLC$. By the main theorem of \cite{W}, the discrete, faithful representations form an open subset of $X(\Si)$. On the other hand, by Mostow rigidity, the discrete faithful representation in $\PSLC$ is unique.  It follows that the character of the representation is isolated in $X(\Si)$, and thus $X^*(W_1) \cap X^*(W_2)$ contains a zero-dimensional component.  Hence (ii) follows from  (i).
\item[(iii)] This is an immediate consequence of Definition 2.1.
\item[(iv)] The quotient map $\SLC \to \PSLC$ induces orientation-preserving maps $X_{\SLC}(\Si) \to X(\Si)$, $X_{\SLC}(W_i) \to X(W_i)$,
and $X_{\SLC}(F) \to X(F)$, where $X_{\SLC}$ denotes the $\SLC$-character variety. Moreover the isotopy $h\colon X^*(F) \to X^*(F)$ used to make $X(W_1)$ transverse to $X(W_2)$ induces an isotopy $h\colon X_{\SLC}^*(F) \to X_{\SLC}^*(F)$ making  $X_{\SLC}(W_1)$ transverse to $X_{\SLC}(W_2)$.
 If $\Si$ is a $\ZZ_2$-homology sphere, then the orientation-preserving map $X_{\SLC}(\Si) \to X(\Si)$ is an isomorphism of character varieties, and the maps
 $X_{\SLC}(W_i) \to X(W_i)$ and $X_{\SLC}(F) \to X(F)$ are local diffeomorphisms in a neighborhood of an isolated intersection point $\chi \in X^*_{\SLC}(W_1) \cap X^*_{\SLC}(W_2)$.  Thus in fact $\la_\PSLC(\Si) = \la_\SLC(\Si).$
\item[(v)] This follows immediately from Definition 2.1.
\item[(vi)] Any representation $\rho:\pi_1(\Si_1 \# \Si_2) \to \PSLC$ is in fact an ordered pair of representations $(\rho_1,\rho_2)$, where
$\rho_i:\pi_1(\Si_i) \to \PSLC$.  Moreover, if $g \in \PSLC$, then $(g \cdot \rho_1,\rho_2)$ is not conjugate to $(\rho_1,\rho_2)$ in general.  In particular, if $\rho_1$ and $\rho_2 $ are both nontrivial, we may use this fact to generate a curve of characters in $X(\Si_1 \# \Si_2)$ containing $\rho$.  Thus, the isolated characters in $X^*(\Si_1 \# \Si_2)$ are characters of representations $\rho = (\rho_{triv},\rho_2)$ with $\rho_2$ an irreducible representation of $\pi_1(\Si_2)$ and characters of representations $\rho = (\rho_1,\rho_{triv})$ with $\rho_1$ an irreducible representation of $\pi_1(\Si_1)$, where $\rho_{triv}$ is the trivial representation.  Hence $$\la_\PSLC(\Si_1 \# \Si_2) =\la_\PSLC(\Si_1) + \la_\PSLC(\Si_2).$$
\end{enumerate}
\end{proof}

We remark that this connected sum formula is more elegant that the corresponding formula for the $\SLC$-invariant. In \cite{BC} it is shown that
$$\la_\SLC(\Si_1 \# \Si_2) =|H^1(\Si_2;\ZZ_2)|\la_\SLC(\Si_1) + |H^1(\Si_1;\ZZ_2)|\la_\SLC(\Si_2).$$
In contrast, the $\PSLC$-invariant is additive.

\section{Computations}

We now study a number of examples which will help clarify the relationship between the $\PSLC$- and $\SLC$-invariants.

In \cite{BC}, the invariant $\la_\SLC(\Si)$ is computed for $\Si$ a Seifert fibered integral homology sphere and for $\Si$ the result of $p/q$-Dehn surgery on a twist knot as long as $p/q$ is not a strict boundary slope for $K_\xi$,  where  $K_\xi$ for $\xi \geq 1$ denotes the $\xi$-twist knot in $S^3$
depicted in Figure 1 below. Moreover, the strict boundary slopes for $K_\xi$ are known and are all even. If $\Si$ is an integral homology sphere, then $\la_\PSLC(\Si) = \la_\SLC(\Si)$ by Theorem 3.1 (iv). Thus we have the following:

\begin{figure}[h]
\begin{center}
\leavevmode\hbox{}
\includegraphics[height=1.3in]{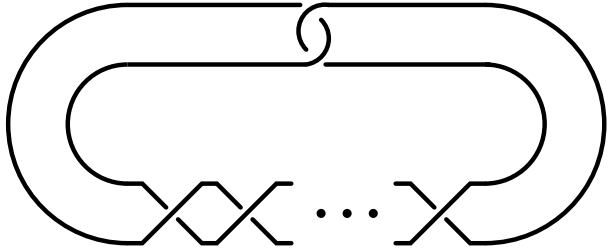}
\caption{{The twist knot $K_\xi$ with a clasp and $\xi$ half
twists.}} \label{ttwist}
\end{center}
\end{figure}

\begin{theorem}
\begin{enumerate}
\item[(i)] Suppose $a_1,\ldots, a_n$ are positive integers that are
pairwise relatively prime, and denote by  $\Si(a_1, \ldots, a_n)$
the associated Seifert fibered homology sphere. Then
$$\la_\PSLC (\Si(a_1,\ldots, a_n)) =\sum_{1 \leq i_1 < i_2 < i_3 \leq n}
 \frac{(a_{i_1}-1)(a_{i_2}-1)(a_{i_3}-1)}{4}.$$
\item[(ii)] If $\Si$ is the result of $p/q$-Dehn surgery on $K_\xi$, where $\xi$ is even and $p$ is odd, then
$$\la_\PSLC(\Si) = \tfrac{1}{4} (\xi |4q+p| + (\xi-2)|p| + 2|2\xi q - p| - 2\xi).  $$
\item[(iii)] If $\Si$ is the result of $p/q$-Dehn surgery on $K_\xi$, where $\xi$ is odd and $p$ is odd, then
$$\la_\PSLC(K_\xi(p/q))= \tfrac{1}{4} \left((\xi-1)(|4q-p| +|p|) + 2|2\xi q+4q-p| -2 \xi\right). $$
\end{enumerate}
\end{theorem}

We now turn our attention to certain Seifert fibered spaces which are not integral homology spheres.  For Seifert fibered spaces which are not $\ZZ_2$-homology spheres, the $\PSLC$- and $\SLC$-invariants may differ. Let $\Si$ be a non-Haken Seifert manifold whose base orbifold is of the form $S^2(p,q,r)$ where $p,q,r \geq 2$.

In \cite{BB}, Proposition D, the authors compute the number of $\PSLC$ characters of $\pi_1\Si$.  They note further that for all irreducible representations
$\rho:\pi_1 \Si \to \PSLC$ the cohomology group $H^1(\Si; Ad \rho)$ is zero. But arguing as in the proof of Theorem 2.1 of \cite{BC}, we find that the character of any such representation in $X^*(W_1) \cap X^*(W_2)$ is an isolated character with multiplicity 1 if it is not the character of a representation which conjugates into $D$, with multiplicity 1/2 if it is the character of a representation which conjugates into $D - Z$, and with multiplicity 1/4 if it is the character of a representation which conjugates into $Z$.  Thus, in this setting, counting the characters of $\pi_1(\Si)$ enables us to compute $\la_\PSLC(\Si)$.

\begin{theorem} Let $\Si$ be a non-Haken Seifert manifold whose base orbifold is of the form $S^2(p,q,r)$ where $p,q,r \geq 2$.  Let $[\cdot]$ denote the greatest integer function.  Then
\begin{eqnarray*}
\la_\PSLC(\Si) &=& \tfrac{1}{4}(p-1)(q-1)(r-1) + \left[\frac{\gcd(p,q)}{2}\right] + \left[\frac{\gcd(p,r)}{2}\right] + \left[\frac{\gcd(q,r)}{2}\right] \\
&& - \left[\frac{\gcd(pq,pr,qr)}{2}\right] - \si(p,q,r),
\end{eqnarray*}
where $\si(p,q,r)$ is 1 if $p,q$, and $r$ are all even and 0 otherwise.
\end{theorem}

\begin{proof}
To see this, we recall Proposition D of \cite{BB}.  It is well-known that $\pi_1(\Si)$ has a presentation of the form
$$\pi_1(\Si) = \langle x, y, h \mid
h \text{ is central, } x^p =h^a, y^q =h^b, (xy)^r =h^c \rangle, $$
where ${\gcd}(a,p) = {\gcd}(b,q) = {\gcd}(c,r)=1$. (See \cite{J} for example.)

\pagebreak
\begin{theorem}[Ben Abdelghani, Boyer, Theorem D, \cite{BB}]
\begin{enumerate}
\item[]
\item[(i)] The number of reducible $\PSLC$-characters of $\pi_1(\Si)$ is
\begin{equation}
\left[ \frac{|aqr+bpr+cpq|}{2}\right] + \left\{ \begin{array}{ll}
                 1 & \mbox{if $\gcd(p,q,r)$ is odd}\\
                 2 & \mbox{if $\gcd(p,q,r)$ is even.}
                 \end{array}\right.
\end{equation}
\item[(ii)] The number of $\PSLC$-characters of $\pi_1(\Si)$ with image a dihedral group of order at least 4 is
\begin{equation}
\si(q,r)\left[\frac{p}{2}\right] + \si(p,r)\left[\frac{q}{2}\right] + \si(p,q)\left[\frac{r}{2}\right] - 2\si(p,q,r)
\end{equation}
where $\si(m,n)$ is one if $m$ and $n$ are both even and zero otherwise.
\item[(iii)] The number of $\PSLC$-characters of $\pi_1(\Si)$ is
\begin{equation}
\begin{split}
\left[\frac{p}{2}\right]\left[\frac{q}{2}\right]\left[\frac{r}{2}\right] &+
\left[\frac{p-1}{2}\right] \left[\frac{q-1}{2}\right] \left[\frac{r-1}{2}\right]
+ \left[ \frac{|aqr+bpr+cpq|}{2}\right] \\
&- \left[ \frac{\gcd(pq,pr,qr)}{2}\right] +\left[\frac{\gcd(p,q)}{2}\right] + \left[\frac{\gcd(p,r)}{2}\right] + \left[\frac{\gcd(q,r)}{2}\right] + 1
\end{split}
\end{equation}
\end{enumerate}
\end{theorem}

Now the dihedral representations are precisely those which conjugate into $D$. Such a representation $\rho$ contributes 1/2
to $\la_\PSLC(\Si)$ if it conjugates into $D - Z$ and contributes 1/4 to $\la_\PSLC(\Si)$ if it conjugates into $Z$.  The reducible characters contribute nothing to $\la_\PSLC(\Si)$, and the remaining characters contribute 1 to $\la_\PSLC(\Si)$.

Note that $$H_1(\Si) \cong \ZZ/m_1 \oplus \ZZ /m_2,$$
where $m_1= \gcd(p, q, r) $  and $m_2=|aqr + bpr + cpq|/m_1.$
Then  $\pi_1(\Si)$ has exactly one character corresponding to a representation which conjugates into $Z$ if $p, q,$ and $r$ are all even and has no such characters if at least one of $p, q,$ and $r$ is odd.  It follows that
$$\la_\PSLC(\Si) = (4.3) - (4.1) - (4.2)/2 - \si(p,q,r)/4.$$
Simplifying this expression, we obtain
\begin{equation*}
\begin{split}
\la_\PSLC(\Si) = \tfrac{1}{4}(p-1)(q-1)(r-1) + \left[\frac{\gcd(p,q)}{2}\right] + \left[\frac{\gcd(p,r)}{2}\right] + \\
 \left[\frac{\gcd(q,r)}{2}\right] - \left[\frac{\gcd(pq,pr,qr)}{2}\right] - \si(p,q,r),
\end{split}
\end{equation*}
as desired.
\end{proof}

We would like to determine $\la_\SLC(\Si)$ for comparison.  As in the $\PSLC$-character variety, every character in $X_{\SLC}^*(W_1) \cap X_{\SLC}^*(W_2)$ is an isolated character, and for any such  character $\chi_{\rho}$ the cohomology group $H^1(\Si; Ad \rho)$ is zero.  Since $\PSLC$ acts freely on the orbit of every irreducible representation, each such character has multiplicity 1 in the intersection cycle $X_{\SLC}^*(W_1) \cdot X_{\SLC}^*(W_2)$.  Thus, $\la_\SLC(\Si)$ is a precise count of the irreducible characters of $\Si$.
  We show
\begin{theorem} Let $\Si$ be a non-Haken Seifert manifold whose base orbifold is of the form $S^2(p,q,r)$ where $p,q,r \geq 2$.  Write $p = 2^{\alpha}p'$, $q = 2^{\beta}q'$, and $r=2^{\gamma}r'$, where $p'$, $q'$, and $r'$ are odd and $\alpha$, $\beta$, and $\gamma$ are greater than or equal to 0. Assume (by reordering  $p$, $q$, and $r$ as necessary) that $\alpha \geq \beta \geq \gamma$. Then
\begin{equation*}
\begin{split}
\la_\SLC(\Si) = \tfrac{1}{4}(p-1)(q-1)(r-1) + \tfrac{\sigma(p,q)}{4}(r-1) + \tfrac{\sigma(p,r)}{4}(q-1) + \tfrac{\sigma(q,r)}{4}(p-1) + \\
\xi_1 \left[\tfrac{\gcd(p,q)}{2}\right] + \xi_2 \left[\tfrac{\gcd(p,r)}{2}\right] +  \xi_2 \left[\tfrac{\gcd(q,r)}{2}\right] - \xi_3 \left[\tfrac{\gcd(pq,pr,qr)}{2}\right] - 4\si(p,q,r),
\end{split}
\end{equation*}
where
$
\begin{cases} \xi_1= 2  & \text{if $\beta > 0$,} \\
\xi_1 = 1 &  \text{otherwise;} \\
\xi_2 = 2  & \text{if  $\gamma>0$ or if $\gamma = 0$ and $\alpha = \beta >0$,} \\
\xi_2 = 1 &  \text{otherwise;}  \\
\xi_3 = 2 & \text{if $\alpha = \beta > \gamma$ or if $\alpha > \beta = \gamma$,} \\
\xi_3 = 1 & \text{otherwise.}
\end{cases}$
\end{theorem}

 In fact, one can check that this equality holds whenever at most one of $p$, $q$, and $r$ is even, regardless of the parities of $a$, $b$, and $c$.

\begin{proof}
Note that if $|aqr+bpr+cpq|$ is odd, then $\Si$ is an $\ZZ_2$-homology sphere, and $\la_\SLC(\Si) = \la_\PSLC(\Si)$ by Theorem 3.1 (iv).  The assertion follows in this case since the formula of Theorem 4.2 agrees with the formula given here if $\beta =0$. We are left with three cases:
\begin{description}
\item[Case 1] $\alpha = \beta = \gamma = 0$; $a$, $b$, and $c$ even
\item[Case 2] $\beta > 0$, $\gamma = 0$
\item[Case 3] $\gamma >0$
\end{description}
We prove Case 3. The remaining cases can be proved using similar arguments; we omit the proofs.  Note that in Case 3, $p$, $q$, and $r$ are even, so $a$, $b$, and $c$ are odd.

We argue along the lines of \cite{BB}.  Note that $h$ is central in $\pi_1(\Si)$.  It is easy to check that if $\rho(h) \neq \pm I$, then $\rho(\pi_1(\Si))$ is diagonalizable.  Thus for any irreducible representation $\rho$ we have $\rho(h) = \pm I$.  We begin by counting the characters of representations with $\rho(h) = I$.

Since $\rho(h) = I$, the irreducible characters are in one-to-one correspondence with the irreducible characters of the triangle group
$$\bigtriangleup (p,q,r) = \langle x,y| x^p = y^q = (xy)^r = 1 \rangle.$$  Note that the character of any representation $\rho$ of $\bigtriangleup (p,q,r)$ with $\rho(x) = \pm I$, $\rho(y) =\pm I$, or $\rho(xy) = \pm I$ is reducible.  Hence the irreducible characters of $\bigtriangleup(p,q,r)$ all lie on the curves $C(j,k)$, $1 \leq j \leq (p/2) - 1$, $1 \leq k \leq (q/2) - 1$, in the character variety $X_{\SLC}(\ZZ_p * \ZZ_q)$, given as follows:  a character $ \chi_z \in C(j,k)$  is the character of a representation $\rho_z$ with
$$
\rho_z(x)  =   \left( \begin{array}{cc}
    \la & 0 \\
    0 & \la^{-1}
    \end{array}\right)
\mbox{and }
\rho_z(y) =  \left( \begin{array}{cc}
    z & 1 \\
    z(\tau - z) - 1 & \tau - z
    \end{array}\right).$$
Here, $\la = e^{2 \pi i j/p}$, $\mu = e^{2 \pi i k/q}$, and $\tau = \mu + \mu^{-1}$.

Each such curve is isomorphic to $\CC$ via the map $z \to \chi_z$.  Moreover the traces $tr(\rho_z(x))$ and $tr(\rho_z(y))$ are constant along the curves, and each curve corresponds to an allowable choice of $tr(\rho(x))$ and $tr(\rho(y))$ for $\rho \in R_{\SLC}(\bigtriangleup(p,q,r))$.

Now $tr(\rho_z(xy)) = z(\la - \la^{-1}) + \la^{-1} \tau$.  Since $\la \neq \la^{-1}$, all possible values of $tr(\rho_z(xy))$ are attained exactly once on the curve - i.e. the map $z \to tr(\rho_z(xy))$ is an isomorphism $\CC \to \CC$.

Note that any irreducible representation $\rho$ of $\bigtriangleup (p,q,r)$ must satisfy $\rho(xy)^r = I$, so $tr(\rho(xy)) = \nu + \nu^{-1}$ for $\nu = e^{2 \pi i l/r}$ for some $l$ with $1 \leq l \leq (r/2) - 1$.   Hence there are precisely $(r/2) - 1$ characters on each curve $C(j,k)$ which are characters in $X_{\SLC}(\bigtriangleup(p,q,r))$.  Summing over all curves $C(j,k)$, we find that we have counted $((p/2) - 1)((q/2) - 1)((r/2) - 1)$ distinct characters.

Unfortunately some of these characters may be reducible.  To count the reducible characters of $\bigtriangleup(p,q,r)$ lying on the curves $C(j,k), 1\leq j\leq (p/2) - 1, 1\leq k \leq (q/2) - 1$, we count the total number of reducible characters in $X_{\SLC}(\bigtriangleup(p,q,r))$ and subtract the number of reducible characters $\chi_{\rho}$ corresponding to representations $\rho$ with $\rho(x) = \pm I$, $\rho(y) = \pm I$, or $\rho(xy) = \pm I.$

Note that $$H_1(\bigtriangleup(p,q,r)) \cong \ZZ / m_1  \oplus \ZZ /m_2' ,$$
where $m_1 = \gcd(p,q,r)$ and $m_2' =  \gcd(pq,pr,qr)/m_1$.
Then the number of reducible characters of $\bigtriangleup(p,q,r)$ is $2 + \gcd(pq,pr,qr)/2$.

If $\rho$ is a reducible character in $X_{\SLC}(\bigtriangleup(p,q,r))$ with $\rho(x) = I$, we see that $\rho(y) = \rho(xy)$.  Since $\rho(y)^q = \rho(xy)^r = I$, we see that there are $(\gcd(q,r)/2) + 1$ such characters.  Similarly there are $(\gcd(q,r)/2) + 1$ characters $\chi_{\rho}$ corresponding to representations $\rho$ with $\rho(x) = -I$, $(\gcd(p,r)/2) - 1$ characters $\chi_{\rho}$ of representations $\rho$ with $\rho(y) = I$ and $\rho(x) \neq \pm I$, and $(\gcd(p,r)/2) - 1$ characters $\chi_{\rho}$ of representations $\rho$ with $\rho(y) = -I$ and $\rho(x) \neq \pm I$.  Finally, there are $(\gcd(p,q)/2) - 1$ characters $\chi_{\rho}$ of representations $\rho$ with $\rho(xy) = I$ and $\rho(x) = \rho(y)^{-1} \neq \pm I$ and $(\gcd(p,q)/2) - 1$ characters $\chi_{\rho}$ of representations $\rho$ with $\rho(xy) = -I$ and $\rho(x) = \rho(y)^{-1} \neq \pm I$.

Thus, there are $\gcd(p,q) + \gcd(p,r) + \gcd(q,r) - 2$ reducible characters with $\rho(x) = \pm I$, $\rho(y) = \pm I$, or $\rho(xy) = \pm I$. It follows that there are
$$4 + (\gcd(pq,pr,qr)/2)  - \gcd(p,q) - \gcd(p,r) - \gcd(q,r)$$
reducible characters of $\Si$ lying on the curves $C(j,k)$, $1 \leq j \leq (p/2) - 1$, $1 \leq k \leq (q/2) - 1$.  Hence there are
$$((p/2) -1)((q/2) - 1)((r/2) - 1) -(\gcd(pq,pr,qr)/2) + \gcd(p,q) + \gcd(p,r) + \gcd(q,r) - 4$$
irreducible characters $\chi$ in $X_{\SLC}(\Si)$ such that $\chi(h) = 2$.

It remains to count the irreducible characters $\chi$ in $X_{\SLC}(\Si)$ with $\chi(h) = -2$.  Note that $a, b,$ and $c$ are odd since $\gcd(p,a) = \gcd(q,b) = \gcd(r,c) = 1$. Then such characters are in one-to-one correspondence with the characters $\chi_{\rho}$ of the triangle group $\bigtriangleup (2p,2q,2r)$ corresponding to representations $\rho$ which satisfy $$\rho(x)^p = \rho(y)^q = \rho(xy)^r = -I.$$  These characters lie on the curves $C(j,k)= \{\rho_z|z\in \CC\}$, $1 \leq j \leq p/2$, $1\leq k \leq q/2$ in $X_{\SLC}(\ZZ_{2p} * \ZZ_{2q})$, where $\la = e^{ \pi i (2j-1)/p}$, $\mu = e^{\pi i (2k - 1)/q}$, and $\rho_z$ is defined as above.  As above, it is easy to check that there are $r/2$ possible choices for $tr(\rho(xy))$ on each curve $C(j,k)$, so there are $(p/2)(q/2)(r/2)$ characters of $\Si$ on these curves.

Again, some of these characters may be reducible.  One checks that there are $\gcd(pq,pr,qr)/2$ reducible characters on these curves if exactly two of $\lcm(p,q,r)/p$, $\lcm(p,q,r)/q$, and $\lcm(p,q,r)/r$ are odd, so either $\alpha = \beta \neq \gamma$ or $\alpha \neq \beta = \gamma$, and none otherwise.  It follows that there are $(p/2)(q/2)(r/2) - \delta \gcd(pq,pr,qr)/2$ irreducible characters $\chi$ of $\Si$ such that $\chi(h) = -2$, where $\delta = 1$ if exactly two of $\lcm(p,q,r)/p$, $\lcm(p,q,r)/q$, and $\lcm(p,q,r)/r$ are odd, so either $\alpha = \beta \neq \gamma$ or $\alpha \neq \beta = \gamma$, and $\delta = 0$ otherwise.  Adding this to the number of irreducible characters $\chi$ of $\Si$ with $\chi(h) = 2$, we find that the total number of irreducible characters is
$$(p/2)(q/2)(r/2) + ((p/2) - 1)((q/2) - 1)((r/2) - 1) - \tfrac{\xi_3}{2} \gcd(pq,pr,qr) $$
 $$+ \gcd(p,q) + \gcd(p,r) + \gcd(q,r)-4.$$
 Finally, note that
 $$(p/2)(q/2)(r/2) + ((p/2) - 1)((q/2) - 1)((r/2) - 1) = \tfrac{1}{4}(p-1)(q-1)(r-1)$$
 $$ + \tfrac{1}{4}(r-1) + \tfrac{1}{4}(q-1) + \tfrac{1}{4}(p-1).$$
This proves Case 3, since each character in $X^*_{\SLC}(\Si)$ contributes 1 to the $\SLC$-Casson invariant.

\end{proof}

\section{Reflections}

We conclude the paper with some observations, questions, and conjectures. A detailed analysis of these points will be the topic of a future paper.

The following statements are either known or should be straightforward to prove using standard techniques.(See \cite{C2}, \cite{BZ2}, and \cite{BZ3} for similar arguments.)
\begin{enumerate}
\item A representation $\rho:\pi_1 \to \PSLC$ lifts to a representation into $\SLC$ if and only if $w_2(\rho)=0 \in H^2(\Si;\ZZ_2) $.
\item If $w_2(\rho)=0$, then there are exactly $|H^1(\Si;\ZZ_2)|$ distinct lifts of $\rho$.  If one lift is $\tilde{\rho}$, then all lifts are given by $\al \tilde{\rho}$ for $\al \in H^1(\Si;\ZZ_2) = Hom(\pi_1(\Si),\ZZ_2)$, where $\al \tilde{\rho} (\gamma) = \al(\gamma) \tilde{\rho}(\gamma)$ for $\gamma \in \pi_1(\Si)$.
\item  Lifts of conjugate representations are conjugate.
\item Distinct lifts of a single irreducible representation $\rho$ with $w_2(\rho)=0$ are not conjugate if the image of $\rho$ does not conjugate into $D$.  If the image of $\rho$ conjugates into $D - Z$, then the lifts of $\rho$ are conjugate in pairs: i.e. each lift $\tilde{\rho}$ of $\rho$ is conjugate to exactly one other lift $\al \tilde{\rho}$ of $\rho$.  If the image of $\rho$ conjugates into $Z$, then  each lift $\tilde{\rho}$ of $\rho$ is conjugate to exactly three other lifts $\al \tilde{\rho}$ of $\rho$.
\item A representation $\rho:\pi_1 \to \PSLC$ with $w_2(\rho)=0$ is reducible if and only if its lifts are reducible.
\item A representation $\rho:\pi_1 \to \PSLC$ with $w_2(\rho)=0$ is isolated if and only if its lifts are isolated.
\end{enumerate}

As long as all of these statements are correct, it is clear that the $\PSLC$-characters $\chi$ with $w_2(\chi)=0$ should contribute $\frac{1}{|H^1(\Si;\ZZ_2)|}\la_\SLC(\Si)$ to the $\PSLC$ invariant. Thus,
\
\begin{conjecture} For any closed, orientable 3-manifold $\Si$,
 $$\la_\PSLC(\Si) =  \frac{1}{|H^1(\Si;\ZZ_2)|}\la_\SLC(\Si)+ \sum_{\chi} n_{\chi},$$ where the sum is taken over all isolated characters $\chi \in X^*(W_1) \cap X^*(W_2)$ with $w_2(\chi) \neq 0$.
\end{conjecture}
Note that a surgery formula for the $\SLC$-Casson invariant exists for certain surgeries.  (See \cite{C} and \cite {C1}.) Thus, the piece of the $\PSLC$-Casson invariant counting characters with $w_2 = 0$ can be easily computed in many cases.

It is clear from the above discussion that we can define invariants $\la_{\al}(\Si)$ counting the $\PSLC$-characters with $w_2 = \al$ for each $\al$ in $H^2(\Si;\ZZ_2)$.  Then $\la_\PSLC(\Si) = \sum_{\al}\la_{\al}(\Si)$.  With this notation, Conjecture 4.1 is the assertion
$\la_0 = \frac{1}{|H^1(\Si;\ZZ_2)|}\la_\SLC(\Si)$.

A similar study relating the $SO(3)$-Casson invariant to the $SU(2)$-Casson invariant was undertaken in \cite{C2}.  In this setting, it was shown that for each element $\al$ of $H^2(\Si;\ZZ_2)$, the contribution $\la_{SO(3),\al}$ of conjugacy classes of $SO(3)$-representations with $w_2 = \al$ to the $SO(3)$-invariant was $\frac{1}{|H^1(\Si;\ZZ_2)|}\la_{SU(2)}(\Si)$.  Thus, in the $SO(3)$-theory, we find that $\la_{SO(3),\al} = \la_{SO(3),0}$ for each $\al \in H^2(\Si;\ZZ_2)$.  Since for closed 3-manifolds $\Si$ we have $H^2(\Si;\ZZ_2) \cong H^1(\Si;\ZZ_2)$, we see that $\la_{SO(3)}(\Si) = \la_{SU(2)}(\Si)$.

Our investigation of Seifert fibered spaces shows that the relationship between the $\PSLC$- and $\SLC$-theories is fundamentally different from the relationship between the $SO(3)$- and $SU(2)$-theories.  For example, let us consider a Seifert fibered space $\Si$ with base orbifold $S^2(p,q,r)$, where $p = 4$, $q=6$, and $r=8$. Theorem 3.2 tells us that $\la_{\PSLC}(\Si)=101/4.$  However by Theorem 3.3, we have  $\la_{\SLC}(\Si)= 30,$ so $\la_0 = 30/4$.  It follows that $\sum_{\al \neq 0} \la_{\al} = 71/4$, so clearly there exists an $\al \in H^2(\Si;\ZZ_2)$ such that $\la_{\al}(\Si) \neq \la_0(\Si)$.  Equivalently, $\la_{\PSLC}(\Si) \neq \la_{\SLC}(\Si).$  This dichotomy between the compact and noncompact theories is most interesting.  We are left with the following problem:

\begin{problem}
Find a surgery formula for $\la_{\al}(\Si)$  for each $\al \in
H^2(\Si;\ZZ_2)$.
\end{problem}

These issues will be addressed in a future paper.

\bigskip \noindent
{\bf Acknowledgements.}
I would like to thank Hans Boden and Andrew Clifford for their useful suggestions.

\end{document}